\documentclass[12pt]{article}
\usepackage{latexsym}
\title{Computing the distribution of the maximum in balls-and-boxes problems,
 with application to clusters of disease cases}
\author{ Warren J. Ewens\\
Department of Biology, University of Pennsylvania\\
Philadelphia, PA 19104-6018\\
{\small\tt <wewens@sas.upenn.edu>} \and Herbert S. Wilf\\Department
of Mathematics, University of Pennsylvania\\Philadelphia, PA
19104-6395\\
{\small\tt <wilf@math.upenn.edu>}}
\newcommand{\eqdef}{\, =\kern -12.7pt\raise 6pt\hbox{{\tiny\textrm{def}}}\,\,}
\begin{document}
\maketitle
\begin{abstract}
We present a rapid method for the exact calculation of the
cumulative distribution function of the maximum of multinomially
distributed random variables. The method runs in time $O(mn)$, where
$m$ is the desired maximum and $n$ is the number of variables. We
apply the method to the analysis of two situations where an apparent
clustering of cases of a disease in some locality has raised the
possibility that the disease might be communicable, and this
possibility has been discussed in the recent literature. We conclude
that one of these clusters may be explained on purely random
grounds, whereas the other may not.
\end{abstract}
\section{Introduction}
It happens, from time to time, that cases of a disease will cluster
both geographically and in time, in a manner which seems not to be
random, and which invites further epidemiological study regarding
communicability of the disease.

Of course mathematics alone cannot answer serious questions of
public health, but it can provide guidelines about what sort of
clustering should be regarded as unusual, and what sort is to be
expected. In particular, the calculation of a $P$-value is required
for an objective assessment of any observed event. In this paper we
provide a rapid and exact $P$-value calculation for the standard
``balls-in-boxes'' model appropriate to the disease clustering
situation.
\section{The model}
Suppose that during a certain time period, a number $r$ of cases
of some disease arise randomly in some large population, such as
that of the the U.S. Let $N$ be size of that population  and
$N_0$ be the population of the community in which the seemingly
large number of cases has occurred.

We think of the entire country as consisting of $n=N/N_0$
identical communities, or cells, each containing $N_0$ people, and
we ask:
\begin{quotation}
\noindent If $r$ cases occur randomly in the populations of $n$
communities of the same size, what is the probability that no
community gets more than $m$ cases of the disease?
\end{quotation}
The standard calculation required to answer this question involves
the ``balls-in boxes'' model, discussed below.  If, for example,
it turns out that it is extremely likely that \textit{some}
community of equivalent size to that where the seemingly large
number of cases occurred, purely by chance, we could conclude that
the observed cluster would not be a cause for further
investigation or suspicion of communicability. Likewise, if it
turns out that it is extremely unlikely that, by chance,
\textit{any} community of the size of that of interest would have
the observed number of cases of the disease, then support would be
given to the possibility of a public health hazard.

\section{The mathematics}
Mathematically speaking, we have $r$ ``balls'' (the disease cases)
being dropped randomly into $n$ labeled ``boxes'' (the communities).
The relevant calculation thus concerns the $P$-value associated with
the box (or boxes) having the largest number of balls in it. It is
well known that the distribution function of the maximum of a number
of random variables changes sharply near the mean of the maximum, so
that an exact rather than an approximate calculation is needed to
find this $P$-value. We provide this exact calculation in this
paper.

The $P$-value associated with an observed value $m$ of cases of
the disease in the community of interest is the probability that
the maximum number of balls in any box in $m$ or more. We find
this probability by first finding the probability that no box
contains more than $m$ balls. Denote this probability by
$P(r,n,m)$.

Now, the probability that there are $r_1$ balls in box 1, and $r_2$
in box 2, and \dots, and $r_n$ in box $n$, is given by the well
known multinomial distribution,
\begin{equation}
\label{eq:multi}
\mathrm{Pr}(r_1,r_2,\dots,r_n)=\frac{1}{n^r}\frac{r!}{r_1!r_2!\dots
r_n!}.\qquad\qquad(r=r_1+\dots+r_n)
\end{equation}
The probability that no box contains more than $m$ balls (i.e., the
cumulative distribution function of the maximum of the $r_i$,
evaluated at $m$) is
\begin{equation}
\label{eq:pdef} P(r,n,m)\eqdef \mathrm{Pr}(\mathrm{all}\ r_i\
\mathrm{are}\ \le m)=\sum_{{0\le r_1,r_2,\dots,r_n\le m}\atop
 {r_1+r_2+\dots+r_n=r}}\frac{1}{n^r}\frac{r!}{r_1!r_2!\dots
r_n!}.
\end{equation}
\section{The computation}
 At first sight the expression (\ref{eq:pdef}) seems appallingly complicated
for exact computation, if $r$ and $n$ are large. Various
approximations, such as the Poisson approximation, have been used by
researchers in order to avoid the apparently tedious computation in
(\ref{eq:pdef}).

However the exact calculation can be completely tamed by two steps.
First we introduce the function
\[e_m(x)=1+x+\frac{x^2}{2!}+\frac{x^3}{3!}+\dots+\frac{x^m}{m!},\]
which is simply the $m$th section of the exponential series.
\textit{Then $P(r,n,m)$ is $r!/n^r$ times the coefficient of $x^r$
in the series $e_m(x)^n$.}
 The question of computing a particular coefficient of a high power
 of a given power series is a well studied problem in
computer science, and the following solution, which makes the
computation quite rapid and easy to program, is taken from \cite{nw}
(chap. 21).

Let $f(x)=\sum_ja_jx^j$ be a given power series and let
$h(x)=f(x)^n$. The question is, if $h(x)=\sum_jh_jx^j$, how can we
economically compute the $h_j$'s from the given $a_j$'s? We begin by
taking logarithms of the equation $h=f^n$, to get
$\log{h(x)}=n\log{f(x)}$. Now differentiate both sides with respect
to $x$ to obtain $h'/h=nf'/f$, and cross multiply to eliminate
fractions, yielding $fh'=nhf'$. Next insert the power series
expansions of the various functions into this equation, and multiply
both sides by $x$, for cosmetic reasons, to get
\[\left(\sum_ja_jx^j\right)\left(\sum_{\ell}\ell h_{\ell}x^{\ell}\right)=
n\left(\sum_jh_jx^j\right)\left(\sum_{\ell}{\ell}a_{\ell}x^{\ell}\right).\]
Finally, equate the coefficients of a given power of $x$, say $x^s$,
on both sides of the last equation, which gives,
\[\sum_{\ell=0}^s\ell h_{\ell}a_{s-\ell}=n\sum_{\ell=0}^s\ell
a_{\ell}h_{s-\ell}.\]

This is a recurrence relation. We can use it to compute the
unknown $h_j$'s successively, in the order $h_0,h_1,h_2,\dots$. To
make this explicit, we can rewrite the above in the form
\begin{equation}
\label{eq:hrec}
h_s=\frac{1}{sa_0}\sum_{\ell=1}^s((n+1)\ell-s)a_{\ell}h_{s-\ell}.\qquad(s=1,2,3,\dots)
\end{equation}
In this form it is clear that each $h_s$ is determined from
$h_0,h_1,\dots,h_{s-1}$.

In the particular case at hand, of powers of the truncated
exponential series $e_k(x)$, we have $a_j=1/j!$, for $0\le j\le m$,
and $a_j=0$ for all other values of $j$. The recurrence takes the
form
\begin{equation}
\label{eq:hmult}h_s=\frac{1}{s}
\sum_{\ell=1}^{\min{(s,m)}}((n+1)\ell-s)\frac{h_{s-\ell}}{\ell
!}.\qquad(s=1,2,3,\dots)
\end{equation}

We summarize the calculation procedure as follows. To compute
$P(r,n,m)$ as defined by eq. (\ref{eq:pdef}) above,
\begin{itemize}
\item Take $h_0=1$ and successively compute $h_1,h_2,\dots,h_r$ from
(\ref{eq:hmult}).
\item Then $P(r,n,m)=r!h_r/n^r$.
\end{itemize}

A remarkable feature of this algorithm is that the computation of
each $h_s$ requires the knowledge of only $m$ earlier values, so the
entire computation can be done with just $m$ units of array storage.
For example, it can find the probability that the maximum is $\le
8$, for 15000 balls in 10000 boxes using only 8 array storage
locations. In summary, it runs in time which is $O(mn)$ and uses
only $O(m)$ storage.

We remark that as we have presented it this method works only for
the situation in which the cells have equal probabilities. It can be
extended, at only a small extra cost, to the case of unequal
probabilities, which may be useful for power calculations.
\section{Some related work}
The problem of finding the distribution of the maximum occupancy in
a balls-and-cells problem is very old. Already in Barton and David
\cite{bd} one finds the first observation above, namely that the
desired probability is a certain coefficient in a power of a given
power series. In \cite{fre} this observation of Barton and David is
cited, and is said to be ``not in a form convenient for computing,''
which is true absent our second step above, in eq. (\ref{eq:hrec})
of vastly accelerating the computation of the high power of the
given series.

Freeman's algorithm in \cite{fre} sought to economize the
computation by grouping together vectors of occupancy numbers which,
as unordered multisets, were the same. Hence he listed partitions
with given largest size part, and counted the occupancies of that
subset of all partitions. This is a large amount more labor than our
method above, which requires computing time roughly proportional to
the square of the number of ``balls,'' whereas earlier methods
required exponential time.

Likewise the recurrence (\ref{eq:hrec}) for computing powers of
power series has a long history. Although we have followed \cite{nw}
in our presentation, the recurrence method was certainly not
invented by them, and is described in several earlier works.
Nonetheless, the concatenation of the two methods in connnection
with finding the distribution of the maximum cell occupancy seems to
be new.

\section{Leukemia: two examples}
{\it Example 1.} We consider first the much discussed case (see
\cite{cwh} and \cite{cwh1}) of  childhood leukemia in Niles, IL in
the five year period 1956--1960. Heath \cite{cwh} gives a total of
eight cases in this town during this period, as compared to an
expected number of 1.6. In 1960 the population of Niles was about
20,000 people. The total population of the U.S. in 1960 was
approximately 180,000,000 people. Therefore the U.S. population in
1960 can be thought of as consisting of 9,000 cells, the population
of each being 20,000 people. An expected number of 1.6 in Niles
would then correspond to a total of about 14,400 cases  in the U.S.
in the five year period studied.

Using the formula above, we therefore computed the exact probability
that if 14,400 balls are distributed randomly into 9000 cells, then
no cell will get more than $m$ balls, for each $m=6,\dots,12$, and
in particular for $m=8$. We also computed the $P$-value for each of
these values of $m$ using the fact that the $P$-values corresponding
to an observed maximum of $m$ is given by $1 - P(14400,9000,m-1)$.

For comparison, we ran a Monte Carlo computer experiment in which we
repeated 1000 times the operation of distributing 14400 balls
randomly into 9000 cells, and recorded the frequencies of the
maximum occupancy numbers, thus giving an empirical distribution
function for $m$. (1000 replications are needed to give an estimate
of the $P$-value for $m=8$ that is accurate to within $\pm 0.01$
with probability 0.95.)  The results of the this simulation, and the
exact $P(14400,9000,m)$ computations are shown below, together with
the exact $P$-values.

\medskip

\begin{center}
\( \begin{array}{lccc} m&P(14400,9000,m)&\mathrm{Monte\ Carlo}&P\mathrm{-value}\\
\hline
6&0.000005&0.000&1.000000\\
7&0.095395&0.096&0.999995\\
8&0.664954&0.678&0.904605\\
9&0.937864&0.944&0.335046\\
10&0.990843&0.993&0.062136\\
11&0.998788&0.998&0.009157\\
12&0.999852&0.999&0.001212
\end{array} \)
\end{center}
The computation of the above table of exact values, on a PC
running the computer algebra system Maple,\footnote{Program
available on request} required less than five seconds. The Monte
Carlo computation required about thirty minutes. In both the Monte
Carlo simulation and the exact calculations we observe the
expected rapid change of $P$-values as $m$ increases, emphasizing
the need for exact $P$-value calculations as discussed above.

We conclude from the above that the probability that some cell of
population 20,000 would have gotten 8 or more cases in the five
year period studied is about 90 percent. Thus the Niles data do
not appear, so far as formal $P$-value calculations are concerned,
to
show a significant cluster of cases of childhood leukemia.\\

\noindent{\it Example 2.} Twelve cases of acute lymphocytic leukemia
were observed \cite{cc} in Churchill County, NV, among persons who
had been residents of the county at the time of diagnosis, in the
three year period 1999-2001. Concern was expressed that this was due
to the exposure to some agent associated with a nearby naval air
station. At that time the county had a population of approximately
24,000. The entire U.S. had a population of approximately
288,000,000, equivalent to 12,000 units, or cells, each of the size
of Churchill County. The State Epidemiologist, Dr. Randall Todd,
estimated that, based on its population, about one case would be
expected in Churchill County every five years. If we use that
estimate, the incidence in the U.S. as a whole would be 12000 cases
per five years, or 8000 cases per three year period.

In this case we need the distribution function of the maximum number
of balls in any cell if 8000 balls are thrown at random into 12,000
cells. The results are as shown below.

\medskip

\begin{center}
\( \begin{array}{lccc} m&P(8000,12000,m)&P\mathrm{-value}\\
\hline
4&0.000472&1.000000\\
5&0.436361&0.999528\\
6&0.925122&0.563639\\
7&0.993604&0.074878\\
8&0.999528&0.006396\\
\end{array} \)
\end{center}

Clearly the observed incidence of twelve cases in Churchill county
cannot reasonably be ascribed to chance, and further epidemiological
investigation is warranted.
\section{Further comments on the $P$-value}
The $P$-value corresponding to any value of $m$ in the
balls-in-boxes case can in principle be calculated exactly using
standard inclusion/exclusion formul\ae. In practice this seems
extremely difficult, because the alternating signs can cause
catastrophic loss of significant digits. A Poisson approximation is
also possible but may be inaccurate, particularly around the tails
of the distribution. Our exact method, described in eq.
(\ref{eq:hmult}) above, is fast and does not suffer from any of
those problems.

A further comment about $P$-values is more wide-ranging. Many
diseases might come to our attention because of an apparent
clustering in some location in some time period. Also, many
different time periods might be potentially observed. An overall
$P$-value calculation, taking these matters into consideration,
would be desirable, but in practice would probably be impossibly
difficult, since no precise value can be attached to ``the number of
diseases that might come to our attention'' or, possibly, to the
number of time periods that we might have considered.
\section{A disclaimer}
Mathematics cannot prove or disprove the communicability of a
disease process. It can only help to define the word ``unusual.''
The benchmark given above seems like an appropriate one to use when
investigating an outbreak which is localized spatially, temporally,
or both. By this benchmark, the clustering of leukemia cases in
Niles, Illinois, between 1956 and 1960 was not unusual. In fact some
collection of that number of cases in some community the size of
Niles, in a five year period of keeping records, was to be expected
with high probability. On the other hand, the Churchill County data
seem extremely significant.

\end{document}